\documentclass[12pt]{article}
\usepackage{amsmath}
\usepackage{amsfonts}
\usepackage{mathrsfs}
\usepackage{amssymb}
\usepackage[usenames]{color}
\usepackage{CJK}

\def\vspb{\vs{5pt}}

\def\qed{\hfill\mbox{$\Box$}}
\def\bs{\backslash}
\def\cl{\centerline}

\def\ni{\noindent}

\def\ssc{\scriptscriptstyle}

\def\vs{\vspace*}
\def\wh{\widehat}

\numberwithin{equation}{section}
\newtheorem{theo}{Theorem}[section]

\newtheorem{coro}[theo]{Corollary}
\newtheorem{lemm}[theo]{Lemma}
\newtheorem{prop}[theo]{Proposition}

\def\a{\alpha}
\def\b{\beta}
\def\d{\delta}

\def\e{\eta}
\def\g{\gamma}

\def\l{\lambda}

\def\s{\sigma}
\def\t{\theta}

\def\C{\mathbb{C}}
\def\Z{\mathbb{Z}}

\def\adb{{\rm ad\ssc\,}{\cal B}(q)}
\def\Autb{{\rm Aut\ssc\,}{\cal B}(q)}
\def\BB{{\cal B}(q)}
\def\BBB{{\cal B}}
\def\Derb{{\rm Der\ssc\,}{\cal B}(q)}
\def\Intb{{\rm Int\ssc\,}{\cal B}(q)}
\def\Vir{\rm Vir}

\begin{document}
\cl{\large\bf {Structure of a class of Lie algebras of Block type
}\footnote{Supported by NSF grant 10825101 of China\\\indent\
$^*$Corresponding author: chgxia@mail.ustc.edu.cn}} \vskip5pt

\cl{{Chunguang Xia$^*$, Taijie You$^\dag$, Liji Zhou$^\dag$}}

\cl{\small \it $^*$Wu Wen-Tsun Key Laboratory of Mathematics,
University of Science and Technology of China\vs{-4pt}}

\cl{\small\it Hefei 230026, China}

\cl{\small \it $^\dag$School of Mathematics and Computer Science,
Guizhou Normal University\vs{-4pt}}

\cl{\small\it Guizhou 550001, China}

\cl{\small\it Email:\ chgxia@mail.ustc.edu.cn, youtj@gznu.edu.cn,
zhou\_li\_ji@126.com}\vs{5pt}

\par\ni {\small{\bf Abstract.} Let $\BB$ be a class of Lie algebras of Block
type with basis $\{L_{\a,i}\,|\,\a,i\in\Z, i\geq 0\}$ and relations
$[L_{\a,i},L_{\b,j}]=\left(\b(i+q)-\a(j+q)\right)L_{\a+\b,i+j}$,
where $q$ is a positive integer. In this paper, it is shown that
$\BB$ are different from each other for distinct positive integers
$q$'s. The automorphism groups, the derivation algebras and the
central extensions of all $\BB$ are also uniformly and explicitly
described, which generalize some previous results. \vskip5pt

\ni{\bf Key words:} Lie algebras of Block type; automorphism;
derivation; central extension.

\ni{\it Mathematics Subject Classification (2000):} 17B05; 17B40;
17B56; 17B65; 17B68.}

\vskip10pt \ni {\bf 1. \
Introduction}\setcounter{section}{1}\setcounter{theo}{0}

\ni In [\ref{B1}], Block introduced a class of infinite dimensional
simple Lie algebras. Partially due to their relations to the
(centerless) Virasoro algebra, generalizations of Lie algebras of
this type (usually referred to as {\it Lie algebras of Block type})
have received many authors' interests (see, e.g., [\ref{DZ},
\ref{OZ}--\ref{SX}, \ref{SZho}, \ref{X1}, \ref{X2}, \ref{Z},
\ref{ZM}]). These Lie algebras are constructed from pairs
$(\mathscr{A},\mathscr{D})$ consisting of a commutative associative
algebra $\mathscr{A}$ with an identity element and a finite
dimensional abelian derivation subalgebra $\mathscr{D}$ such that
$\mathscr{A}$ is $\mathscr{D}$-simple (such pairs are classified in
[\ref{SXZ}]). The representation theory for the simple Lie algebras
of Block type is far from being well developed, except for some
quasifinite representations (see, for example, [\ref{S1}--\ref{S3},
\ref{WT1}]), which is partially because the structure theory for the
Lie algebras of Block type is not well developed yet.

In the present paper, we concentrate on a class of Lie algebras
$\BB$, where $q$ is a positive integer, with basis
$\{L_{\a,i}\,|\,\a\in\Z,i\in\Z_+\}$ and brackets
\begin{equation}\label{B-block}
[L_{\a,i},L_{\b,j}]=\left(\b(i+q)-\a(j+q)\right)L_{\a+\b,i+j} \mbox{
\ for \ } \a,\b\in\Z,\,\,i,j\in\Z_+.
\end{equation}
The structure theory, including automorphism group, derivation
algebra and central extension, for the special case $\BBB(1)$ were
considered in [\ref{WT2}, \ref{XW}]. Our goal is to study the
structure theory (and later the representation theory) for the whole
class of Lie algebras $\BB$ not only for a particular Lie algebra
$\BBB(1)$. Our first result is the following.

\begin{theo}[Isomorphism Theorem]\label{thm-isomorphism}
Lie algebras $\BB$ are different from each other for distinct
positive integers $q$'s, namely,
$$
\BBB(q_1)\cong\BBB(q_2)\ \ \Longleftrightarrow\ \ q_1=q_2.
$$
\end{theo}

The quasifinite representations of Lie algebras of Block type were
initially studied in [\ref{S1}, \ref{S2}], which seem to be
intrigued by Mathieu's classification of Harish-Chandra modules over
the well-known Virasoro algebra in [\ref{M}]. As pointed in
[\ref{WT1}], the central extensions of $\BBB(1)$, which can be
realized as a subalgebra of the Block type Lie algebra $\BBB$
studied in [\ref{S1}] and which is quite different from  $\BBB$,
contains a subalgebra isomorphic to the Virasoro algebra, thus the
representations for $\BBB(1)$ and its central extensions will be
much more interesting. We will show in Theorem \ref{thm-2cocycle}
that the central extension, denoted $\wh\BB$, of $\BB$ for $q\ge1$
is given by
$$
[L_{\a,i},L_{\b,j}]=\left(\b(i+q)-\a(j+q)\right)L_{\a+\b,i+j}
+\d_{\a+\b,0}\d_{i,0}\d_{j,0}\frac{\a^3-\a}{12}c,
$$
for $\a,\b\in \Z,\,\,i,j\in\Z_+$. One sees that  $\wh\BB$ contains a
subalgebra with basis $\{q^{-1}L_{\a,0}, c\,|\,\a\in\Z\}$ isomorphic
to the Virasoro algebra. One also sees that $\BB$ contains a
subalgebra with basis $\{q^{-1}L_{\a,qi}\,|\,\a\in\Z,i\in\Z_+\}$,
which is isomorphic to $\BBB(1)$. Due to these two facts and the
above isomorphism theorem, one may expect that the representation
theory of $\wh\BB$ ($q\ge2$) may be more interesting than that of
$\wh{\BBB(1)}$. We would also like to point out that although $\BB$
is $\Z$-graded with respect to eigenvalues of ${\rm ad}_{L_{0,0}}$,
it is not finitely-generated $\Z$-graded, thus some classical
methods due to Farnsteiner [\ref{F}] (which are efficient for
finitely generated graded Lie algebras), cannot be applied in our
case here.

Now we outline our main results in the present paper. In Section 2,
after giving the proof of Theorem \ref{thm-isomorphism}, we describe
the automorphism group of $\BB$ (Theorem \ref{thm-automorphiam}),
which in particular shows that $\BB$ has no nontrivial inner
automorphisms. In Section 3, by employing a technique developed in
[\ref{SZho}], we characterize the structure of the derivation
algebra of $\BB$ and prove that the first cohomology group of $\BB$
with coefficients in its adjoint module is one-dimensional (Theorem
\ref{thm-derivation}). Finally in Section 4, we determine the
central extensions and the second cohomology group of $\BB$  and
prove that $\wh\BB$ is an essentially unique nontrivial
one-dimensional central extension of $\BB$ (Theorem
\ref{thm-2cocycle}).

Throughout this paper, $q$ will denote a fixed positive integer. All
the vector spaces are assumed over the complex field $\C$. As usual,
we denote by $\Z$ the ring of integers and by $\Z_+$ the set of
nonnegative integers.

\vskip15pt \ni {\bf 2. \ Automorphisms of
$\BB$}\setcounter{section}{2}\setcounter{theo}{0}\setcounter{equation}{0}

\ni {An element $F\in\BB$ is called
\baselineskip3pt\lineskip7pt\parskip-3pt
\begin{itemize}\parskip-3pt
  \item[\rm (1)] {\it $\rm ad$-locally finite} if for any given $v\in \BB$,
the subspace ${\rm Span}\{{\rm ad}_F^m( v)\,|\,m\in\Z_+\}$ of $\BB$
is finite dimensional;
  \item[\rm (2)] {\it $\rm ad$-locally nilpotent} if for any given $v\in \BB$, there
exists some $N>0$ such that ${\rm ad}_F^N( v)=0$\vspace*{-8pt}.
\end{itemize}}
\ni Denote by $\Autb$ the \emph{automorphism group} of $\BB$, and
$\Intb$ the \emph{inner automorphism group} of $\BB$, i.e., the
subgroup of $\Autb$ generated by ${\rm exp}^{{\rm ad}_x}$ for $\rm
ad$-locally nilpotent elements $x$'s. Recall that the centerless
Virasoro algebra $\Vir$ (or Witt algebra) with basis
$\{L_\a\,|\,\a\in\Z\}$ is defined by the following commutation
relations:
$$
[L_\a,L_\b]=(\b-\a)L_{\a+\b}\ \ \mbox{for}\ \ \a,\b\in\Z.
$$

In this section, we first prove that $\BB$ contains a unique locally
finite element $L_{0,0}$ (up to scalars), which is not locally
nilpotent (thus $\Intb$ is trivial).  Then we introduce a lemma
which gives some relations between $\Vir$ and $\BB$. After that, by
giving the structure of automorphism group of $\Vir$, and
introducing two useful lemmas, we present a  proof of Theorem
\ref{thm-isomorphism}. Finally we completely characterize the
structure of $\Autb$ (Theorem \ref{thm-automorphiam}).

\begin{lemm}\label{lemma-aut-local-finite}
We have the following facts:
\baselineskip3pt\lineskip7pt\parskip-3pt
\begin{itemize}\parskip-3pt
  \item[{\rm(1)}] Up to scalars, $L_{0,0}$ is  the unique locally finite element of
$\BB$.
  \item[{\rm(2)}] $L_{0,0}$ is not locally
nilpotent, thus $\Intb$ is trivial.
\end{itemize}
\end{lemm}

\ni{\it Proof.} \ \ (1) Introduce a {\it  total order} $\prec$ on
$\Z\times\Z_+$ defined by
$$
(\a,i)\prec(\b,j)\ \ \mbox{if}\ \ \a<\b\ \ \mbox{or}\ \ \a=\b, i>j.
$$
Let $F=\sum_{(\a,i)\in I_F}\l_{\a,i}L_{\a,i}$ be any locally finite
element of $\BB$, where $I_F$ is a finite subset of $\Z\times\Z_+$.
First assume that there exists $\l_{\a,i}\neq0$ for some $\a< 0$.
Set
$$
(\a_0,i_0)={\rm min}\{(\a,i)\in I_F\,|\,\l_{\a,i}\ne0\} \mbox{
\ (under the sense of total order $\prec$)}.
$$
Then clearly
$\a_0<0$. By rescaling $F$, we may suppose
$$
F=L_{\a_0,i_0}+\sum_{(\a_0,i_0)\prec(\a,i)}\l_{\a,i}L_{\a,i}.
$$
In this case we say that $F$ has the {\it minimal term}
$L_{\a_0,i_0}$. Recall that $[L_{\a_0,i_0},L_{\b,j}]=f(\b,j)
L_{\a_0+\b,i_0+j}$, where $f(\b,j):=\b(i_0+q)-\a_0(j+q)$. For any
given $j_0\in\Z_+$, we can choose small enough $\b_0$ such that
$f(\b_0,j_0)<0$ and then
$$
f(\b_0+k\a_0,j_0+ki_0)=f(\b_0,j_0)+k\a_0 q<0\,\ \ \mbox{for all}\ \
k\in\Z_+,
$$
which implies that ${\rm ad}_F^k(L_{\b_0,j_0})$, with minimal terms
$L_{\b_0+k\a_0,j_0+ki_0}$, are linear independent for all $k$, thus
$F$ is not $\rm ad$-locally finite. Hence $\l_{\a,i}=0$ for all
$\a<0$. Similarly, we can prove that $\l_{\a,i}=0$ for all $\a>0$.
In turn, we can rewrite $F=\sum_{i\in I'_F}\l_{0,i} L_{0,i}$, where
$I'_F$ is a finite subset of $\Z_+$. If there exists $\l_{0,i}\neq0$
for some $i>0$, then similarly we can assume
$$
F=L_{0,i_0}+\sum_{i_0>i}\l_{0,i} L_{0,i}.
$$
Note that $[L_{0,i_0},L_{\b,j}]=g(\b)L_{\b,i_0+j}$, where
$g(\b):=\b(i_0+q)$. In particular, we  see that $g(1)=i_0+q>0$,
which implies that ${\rm ad}_F^k(L_{1,j})$ are linear independent
for all $k$. So $\l_{0,i}=0$ for all $i>0$. Thus $F=\l_{0,0}
L_{0,0}$ for some $\l_{0,0}\in\C$, namely, $L_{0,0}$ is up to
scalars the unique locally finite element of $\BB$.

(2) Note that any locally nilpotent element must be locally finite
element by definition. Since ${\rm ad}^N_{L_{0,0}}L_{1,0}=q^N
L_{1,0}\ne0$ for any positive integer $N$, we see that  $L_{0,0}$ is
not locally nilpotent. By (1), $\BB$ does not contain any nonzero
locally nilpotent elements (in particular, the inner automorphism
group of $\BB$ is trivial). \qed\vskip5pt

\begin{lemm}\label{lemm-aut-Sub-Vir}
Let ${\rm Span}\{{\bar L}_\a\,|\,\a\in\Z\}$ be a subalgebra of
$\BB$, which is isomorphic to $\Vir$, namely $[{\bar L}_\a,{\bar
L}_\b]=(\b-\a){\bar L}_{\a+\b}$. If ${\bar L}_0\in\C L_{0,0}$, then
there exists some $s_0\in\Z^*=\Z\backslash\{0\}$ such that  ${\bar
L}_\a\in\C L_{s_0\a,0}$ for all $\a\in\Z$.
\end{lemm}

\ni{\it Proof.}\ \ Assume ${\bar L}_0=q^{-1}a_0L_{0,0}$ for some
$a_0\in\C^*$. Let $0\ne\a\in\Z$. Write ${\bar
L}_{\a}=\sum_{(\b,j)\in J_\a}\mu^\a_{\b,j}L_{\b,j}$, where $J_\a$ is
a finite subset of $\Z\times\Z_+$. Then
$$
\a \sum_{(\b,j)\in J_\a}\mu^\a_{\b,j}L_{\b,j}=\a {\bar L}_\a=[{\bar
L}_0,{\bar L}_\a]= \Big[q^{-1}a_0L_{0,0},\sum_{(\b,j)\in
J_\a}\mu^\a_{\b,j}L_{\b,j}\Big] =\sum_{(\b,j)\in J_\a}\b
a_0\mu^\a_{\b,j}L_{\b,j},
$$
which implies
\begin{equation}\label{equ-aut-0a}
(\a-\b a_0)\mu^\a_{\b,j}=0\mbox{ \ for all \ $0\ne\a\in\Z$ and
$(\b,j)\in J_\a$.}
\end{equation}
Obviously $\mu^\a_{0,j}=0$ for all $0\ne\a\in\Z$ and $(0,j)\in
J_\a$. Note that for any nonzero $\a$, there exists at least one
$\mu^\a_{\b,j}\ne0$ with $\b\ne0$. In particular, take $\a=1$, there
exist $s_0\in\Z^*$ and some $j_0$ such that $\mu^1_{s_0,j_0}\ne0$,
and so $a_0=\frac1{s_0}$. Thus \eqref{equ-aut-0a} implies in general
$\mu^\a_{\b,j}=0$ if $\b\ne s_0\a$.
 Hence we can rewrite ${\bar L}_{\a}=\sum_{j\in
J'_\a}\mu^\a_j L_{s_0\a,j}$, where $\mu^\a_j=\mu^\a_{s_0\a,j}$ and
$J'_\a=\{j\,|\,(s_0\a,j)\in J_a\}\subset\Z_+$. Then
\begin{eqnarray}\label{equ-aut--aa}
 2\a (qs_0)^{-1} L_{0,0}&=&[{\bar L}_{-\a},{\bar L}_\a]=\Big[\sum_{i\in J'_{-\a}}\mu^{-\a}_i
L_{-s_0\a,i},\sum_{j\in J'_\a}\mu^\a_j L_{s_0\a,j}\Big]\nonumber\\
&=&\!\!\!\sum_{(i,j)\in J'_{-\a}\times J'_\a}\!\!\!(i+j+2q)\a
s_0\mu^{-\a}_i \mu^\a_j L_{0,i+j}.
\end{eqnarray}
Let $i_0={\rm max}\{i\in J'_{-\a}\,|\,\mu^{-\a}_i\ne0\}$,
$j_{0}={\rm max}\{j\in J'_\a\,|\,\mu^\a_j\ne0\}$. If $i_0+j_0>0$,
then the right hand side of \eqref{equ-aut--aa} contains the nonzero
term $(i_0+j_0+2q)\a s_0\mu^{-\a}_{i_0} \mu^\a_{j_0}
 L_{0,i_0+j_0}$, which is
not in $\C L_{0,0}$ (since $i_0,j_0\ge0$). Thus $i_0=j_0=0$, in
particular, ${\bar L}_\a\in\C L_{s_0\a,0}$. \qed\vskip5pt

\begin{prop}[\mbox{See for example [\ref{PZ}, \ref{SZha}]}]\label{prop-aut-vir}
The following two maps are automorphisms of $\Vir$:
\begin{eqnarray}
&\bar\t_\mu: & {\Vir\rightarrow\Vir}, \ \ \ \ L_\a\mapsto\mu^\a L_\a, \nonumber\\
&\bar\d_s: & {\Vir\rightarrow\Vir}, \ \ \ \ L_\a\mapsto s L_{s\a},
\nonumber
\end{eqnarray}
where $\mu\in\C^*=\C\bs\{0\}$, $s\in\{\pm1\}\cong\Z/2\Z$.
Furthermore, the structure of the automorphism group of $\Vir$,
denoted ${\rm Aut}(\Vir)$, is given by ${\rm
Aut}(\Vir)\cong\C^*\rtimes\Z/2\Z.$
\end{prop}

Let $\BBB(q_1)$ and $\BBB(q_2)$ be two Lie algebras of Block type
defined as \eqref{B-block}.  We use the notations $L_{\a,i}$ and
$L'_{\a,i}$ to stand for the base elements of $\BBB(q_1)$ and
$\BBB(q_2)$ respectively. Suppose
$\tau:\BBB(q_1)\rightarrow\BBB(q_2)$ \mbox{is a Lie algebra
isomorphism}. We want to prove $q_1=q_2$. First we give the
following two useful lemmas.

\begin{lemm}\label{lemma-aut-a0}
We have $\tau(L_{\a,0})=sq_1q_2^{-1}\mu^\a L'_{s\a,0}$, where
$\mu\in\C^*$, $s\in\{\pm1\}$.
\end{lemm}

\ni{\it Proof.} \ \  Let $L'_\a=\tau(q_1^{-1}L_{\a,0})$ for
$\a\in\Z$. Since ${\cal V}_1={\rm
Span}\{q_1^{-1}L_{\a,0}\,|\,\a\in\Z\}$ is a subalgebra of
$\BBB(q_1)$ isomorphic to $\Vir$, we see that ${\cal V}_2={\rm
Span}\{L'_\a\,|\,\a\in\Z\}$ is a subalgebra of $\BBB(q_2)$
isomorphic to $\Vir$. Furthermore, since $L_{0,0}$ and $L'_{0,0}$
are up to scalars the unique ${\rm ad}$-locally finite elements in
$\BBB(q_1)$ and $\BBB(q_2)$ respectively, we must have
$L'_0=\tau(q_1^{-1}L_{0,0})\in\C L'_{0,0}$. Hence Lemma
\ref{lemm-aut-Sub-Vir} implies that there exists $s\in\Z^*$ such
that $L'_\a\in\C L'_{s\a,0}$. Analogously, there exists $s'\in\Z^*$
such that $\tau^{-1}(L'_{\a,0})\in\C L_{s'\a,0}$ for all $\a\in\Z$.
In particular $ss'=1$. Thus $s=\pm1$.

We may assume that $s=1$ as the case for $s=-1$ is similar. Thus we
can suppose $\tau(L_{\a,0})=q_1\mu_\a L'_{\a,0}$ for some
$\mu_\a\in\C^*$ for $\a\in\Z$. Applying $\tau$ to
$[L_{1,0},L_{\a,0}]=(\a-1)q_1 L_{\a+1,0}$ gives
$(\a-1)(\mu_{\a+1}-q_2\mu_1\mu_\a)=0$, which by induction implies
\begin{equation}\label{equ-aut-mua}
\mu_\a=\left\{
\begin{aligned}
\mu_2(q_2\mu_1)^{\a-2}\ \ \mbox{if}\ \ \a\ge2,\\
q_2^{-1}(q_2\mu_1)^\a\ \ \mbox{if}\ \ \a\le1.
\end{aligned}\right.
\end{equation}
Applying $\tau$ to $[L_{-2,0},[L_{2,0},L_{3,0}]]=7q_1^2 L_{3,0}$, we
obtain $\mu_{-2}\mu_2 q_2^2=1$, which together with
\eqref{equ-aut-mua} gives $\mu_\a=q_2^{-1}(q_2\mu_1)^\a$. Hence
$\tau(L_{\a,0})=q_1q_2^{-1}(q_2\mu_1)^\a L'_{\a,0}$ for $\a\in\Z$.
\qed

\begin{lemm}\label{lemma-aut-0i}
We have $q_1|q_2$ and $\tau(L_{0,i})=\nu_i L'_{0,q_1^{-1}q_2 i}$,
where $\nu_i\in\C^*$.
\end{lemm}

\ni{\it Proof.} \ \ Assume that $\tau(L_{0,i})=\sum_{(\b,j)\in
K_i}\nu^{\,i}_{\b,j}L'_{\b,j}$, where $K_i$ is some finite subset of
$\Z\times\Z_+$. Applying $\tau$ to the equation
$[L_{0,0},L_{0,i}]=0$, we obtain by Lemma \ref{lemma-aut-a0} that
$$
\sum_{(\b,j)\in K_i}s q_1\b\nu^i_{\b,j}L'_{\b,j}=0,
$$
which implies that $\nu^{\,i}_{\b,j}=0$ if $\b\ne0$. In turn,
$\tau(L_{0,i})$ can be rewritten as
$$
\tau(L_{0,i})=\sum_{j\in K'_i}\nu^{\,i}_jL'_{0,j},\ \ \mbox{where}\
\ \nu^i_j=\nu^i_{0,j},\mbox{ and }K'_i=\{j\,|\,(0,j)\in K_i\}.
$$
Applying $\tau$ to
$[L_{-1,0},[L_{1,0},L_{0,i}]]=-(i+q_1)(i+2q_1)L_{0,i}$, we obtain by
Lemma \ref{lemma-aut-a0} that
$$
\sum_{j\in K'_i}(q_1 j-q_2 i)(q_1 j+q_2 i+3 q_1
q_2)\nu^{\,i}_jL'_{0,j}=0.
$$
Since $i,j\ge0$ and $q_1,q_2>0$, the above implies that
$\nu^{\,i}_{j}=0$ if $j\ne q_1^{-1}q_2 i$. If $q_1\!\nmid\! q_2$,
then we see that $\tau(L_{0,1})=0$, which is impossible since $\tau$
is an isomorphism. Therefore $q_1|q_2$ and $\tau(L_{0,i})=\nu_i
L'_{0,q_1^{-1}q_2 i}$, where
$\nu_i=\nu^{\,i}_{q_1^{-1}q_2i}\in\C^*$. This completes the proof.
\qed\vspb

\ni{\it Proof of Theorem \ref{thm-isomorphism}.}\ \ Now we give the
proof of Theorem \ref{thm-isomorphism}.  Let $\tau$ be a Lie algebra
isomorphism from $\BBB(q_1)$ to $\BBB(q_2)$. By Lemma
\ref{lemma-aut-0i}, we have $q_1|q_2$. On the other hand,
$\tau^{-1}$ is also an Lie algebra isomorphism from $\BBB(q_2)$ to
$\BBB(q_1)$, then $q_2|q_1$. Since $q_1,q_2>0$, we have
$q_1=q_2$.\qed

\begin{coro}\label{coro-automorphiam}
Let $\s\in\Autb$. There exist $\mu\in\C^*$, $s\in\{\pm1\}$ and
$\nu_i\in\C^*$ for $i\in\Z_+$, such that
\baselineskip3pt\lineskip7pt\parskip-3pt
\begin{itemize}\parskip-3pt
  \item[{\rm(1)}] $\s(L_{\a,0})=s\mu^\a
  L_{s\a,0}$,
  \item[{\rm(2)}] $\s(L_{0,i})=\nu_i
  L_{0,i}$.
\end{itemize}
\end{coro}

\ni{\it Proof.} \ \ Parts (1) and (2) follow directly from Lemma
\ref{lemma-aut-a0} and \ref{lemma-aut-0i} respectively (one will see
a further description of $\s(L_{0,i})$ in the following
theorem).\qed

Motivated by the Proposition \ref{prop-aut-vir} and Corollary
\ref{coro-automorphiam}, for any $\mu,\nu\in\C^*,\,s\in\{\pm1\}$,
one can define the following three kinds of linear transformations
of $\BB$:
\begin{eqnarray}
&\t_\mu: & \BB\rightarrow\BB \ \ \ \ L_{\a,i}\mapsto \mu^\a L_{\a,i}; \nonumber\\
&\e_\nu: & \BB\rightarrow\BB \ \ \ \ L_{\a,i}\mapsto \nu^i L_{\a,i}; \nonumber\\
&\d_s: & \BB\rightarrow\BB \ \ \ \ L_{\a,i}\mapsto s
L_{s\a,i}.\nonumber
\end{eqnarray}
One can easily check that they are all automorphisms of $\BB$.
Furthermore, we have the following facts.
{\baselineskip3pt\lineskip7pt\parskip-3pt
\begin{itemize}\parskip-3pt
  \item $\{\t_\mu\,|\,\mu\in\C^*\}\cong\C^*$ is a subgroup of
$\Autb$, where $\t_{\mu_1}\t_{\mu_2}=\t_{\mu_1\mu_2}$ for
$\mu_1,\mu_2\in\C^*$.
  \item $\{\e_\nu\,|\,\nu\in\C^*\}\cong\C^*$ is a subgroup of
$\Autb$, where $\e_{\nu_1}\e_{\nu_2}=\e_{\nu_1\nu_2}$ for
$\nu_1,\nu_2\in\C^*$.
  \item $\{\d_s\,|\,s=\pm1\}\cong\Z/2\Z$ is a subgroup of $\Autb$.
\end{itemize}}

\begin{theo}\label{thm-automorphiam}
Let $\s\in\Autb$, we have $\s(L_{\a,i})=s \mu^\a \nu^i L_{s\a,i}$
for some $\mu,\nu\in\C^*$, $s\in\{\pm1\}$. In particular, $
\Autb\cong(\C^*\times\C^*)\rtimes\Z/2\Z$.
\end{theo}

\ni{\it Proof.} \ \ Let $\s\in\Autb$, by Corollary
\ref{coro-automorphiam}, we have $\s(L_{\a,0})=s\mu^\a L_{s\a,0}$
and $\s(L_{0,i})=\nu_i L_{0,i}$ for some $\mu, \nu_i\in\C^*,
s\in\{\pm1\}$. Applying $\s$ to the equation
$[L_{\a,0},L_{0,i}]=-\a(i+q)L_{\a,i}$ gives
$$
\s(L_{\a,i})=\mu^\a \nu_i L_{s\a,i}\ \ \mbox{for}\ \ \a\in\Z,
i\in\Z_+.
$$
Applying $\s$ to $[L_{0,1},L_{1,i}]=(1+q)L_{1,i+1}$, we obtain
$\nu_{i+1}=s\nu_1\nu_i$, which implies that $\nu_i=s\nu^i$ for some
$\nu\in\C^*$. In turn, $\s(L_{\a,i})=s\mu^\a \nu^i L_{s\a,i}$. In
particular, based on facts of the three maps $\t_\mu, \e_\nu$ and
$\d_s$ stated above, we see that
$\Autb\cong(\C^*\times\C^*)\rtimes\Z/2\Z$. \qed

\vskip15pt \ni {\bf 3. \ Derivations of
$\BB$}\setcounter{section}{3}\setcounter{theo}{0}\setcounter{equation}{0}

\ni Note that $\BB=\oplus_{\a\in\Z}\BB_\a$ is a $\Z$-graded Lie
algebra, where $\BB_\a=\mbox{span}\{L_{\a,i}\,|\,i\in\Z_+\}$. Recall
that a \emph{derivation} $d$ of $\BB$ is a linear transformation on
$\BB$ such that
$$
d([x,y])=[d(x),y]+[x,d(y)] \mbox{ \ for \ } x,y\in\BB.
$$
Denote $\Derb$ and $\adb$ the space of the derivations and
\emph{inner derivations} of $\BB$, respectively. Elements in
$\Derb\bs\adb$ are called \emph{outer derivations}. We say that a
derivation $d$  has \emph{degree} $\a$ $({\rm deg}(d)=\a)$  if
$d\ne0$ and $d(\BB_\b)\subset\BB_{\a+\b}$ for any $\b\in\Z$. Note
that
$$
H^1(\BB)=\Derb /\adb
$$
is the \emph{first cohomology group of $\BB$ with coefficients in
its adjoint module}.

To begin with, we introduce the following notations ($\a\in\Z,
i\in\Z_+$)
\begin{align*}
&\BB^{[j]}=\sum_{\a\in\Z}\BB_\a^{[j]}\ \ \mbox{with}\ \
 \BB_\a^{[j]}=\mbox{span}\{L_{\a,i}\,|\,i\leq j\},\\
&\BB^{(j)}=\sum_{\a\in\Z}\BB_\a^{(j)}\ \ \mbox{with}\ \
 \BB_\a^{(j)}=\mbox{span}\{L_{\a,i}\,|\,i<j\},\\
&\BB_\a=\sum_{j\in\Z_+}\BB_\a^{[j]},\ \ \
(\Derb)_\a=\{d\in\Derb\,|\,{\rm deg}(d)=\a\}.
\end{align*}
Clearly, we have a derivation $D_0$ with degree zero of $\BB$
defined by
\begin{equation}\label{equ-der-outer}
D_0:L_{\b,j}\mapsto j L_{\b,j} \mbox{ \ for \ } \b\in\Z, j\in\Z_+,
\end{equation}
which can be easily checked to be an outer derivation. The main
results of this section is given as follows.

\begin{theo}\label{thm-derivation}
We have $\Derb=\adb \oplus {\cal D}$, where ${\cal D}=\C D_0$. In
particular, the first cohomology group of $\BB$ is one-dimensional,
namely ${\rm dim\ssc\,}H^1(\BB)=1$.
\end{theo}

\ni{\it Proof.} \ \ Let $d\in\Derb$. The proof of the theorem is
equivalent to proving that $d$ is spanned by an inner derivation
ad$_u\in\adb$ for some $u\in\BB$ and $D_0\in {\cal D}$. This will be
done by Lemmas \ref{lemma-der-0st}--\ref{lemma-der-4st}.\qed
\vskip5pt

For a fixed integer $\a\in\Z$, consider a nonzero derivation
$d\in(\Derb)_\a$ such that
\begin{equation}\label{equ-der-suppose}
d(\BB^{[j]})\subset\BB^{[i+j]}\ \ \mbox{for any}\ \ j\in\Z_+,
\end{equation}
where $i\in\Z$ is assumed to be the  minimal integer
satisfying (\ref{equ-der-suppose}). Then we can write
\begin{equation}\label{equ-der-suppose+}
d(L_{\b,j})\equiv e_{\b,j}L_{\a+\b,i+j}(\mbox{mod }\BB^{(i+j)}),
\end{equation}
where $e_{\b,j}\in\C$ and we adopt the convention that if a notation
is not defined but technically appears in an expression, we always
treat it as zero; for example, if $i<0$ in (\ref{equ-der-suppose}),
then $e_{\b,0}=0$ for any $\b\in\Z$.

\begin{lemm}\label{lemma-der-0st}
The minimal integer $i$ satisfying {\rm(\ref{equ-der-suppose})} must
be nonnegative.
\end{lemm}

\ni{\it Proof.} \ \ If not so, then $i<0$. Then
\eqref{equ-der-suppose+} in particular implies $d(L_{\g,0})=0$ for
all $\g\in\Z$. Applying $d$ to $[L_{\b,j},L_{\g,0}]=\left(\g(j+q)-\b
q\right)L_{\b+\g,j}$, we obtain
\begin{equation}\label{equ-der-nonnegtivity}
(\g(i+j+q)-(\a+\b)q)e_{\b,j}=(\g(j+q)-\b q)e_{\b+\g,j}.
\end{equation}
Taking $\g=0$ gives $\a e_{\b,j}=0$. Using this and replacing
$(\b,\g)$ by $(\b,1)$ and $(\b+1,-1)$ respectively in
\eqref{equ-der-nonnegtivity}, we obtain two equations on $e_{\b,j}$
and $e_{\b+1,j}$. Solving these two equations, we get
\begin{equation}\label{equ-der-0-ebj}
(i+2j+3q)e_{\b,j}=0.
\end{equation}
Furthermore, replacing $(\b,\g)$ by $(\b,1)$, $(\b+1,1)$ and
$(\b,2)$ respectively in \eqref{equ-der-nonnegtivity}, we obtain
three equations on $e_{\b,j}$, $e_{\b+1,j}$ and $e_{\b+2,j}$.
Solving the three equations, together with \eqref{equ-der-0-ebj},
gives $2(j+q)(j+2q)(2j+3q)e_{\b,j}=0$, which implies that $e_{\b,j}
=0$ for all $\b\in\Z, j\in\Z_+$, contradicting the minimality of $i$
in (\ref{equ-der-suppose}).\qed

\begin{lemm}\label{lemma-der-1st}
If $\a\neq0$ or $\a=0$, $i\neq0$, then $d$ in
$(\ref{equ-der-suppose})$ is an inner derivation.
\end{lemm}

\ni{\it Proof.} \ \ Applying $d$ to $[L_{\b,j},L_{0,0}]=-\b q
L_{\b,j}$, we have
\begin{equation}\label{equ-der-1-ebj}
\a q e_{\b,j} =(\a(j+q)-\b(i+q))e_{0,0}.
\end{equation}
First suppose $\a\neq0$. Set $u_1=(\a q)^{-1}e_{0,0}L_{\a,i}\in\BB$
and let $d'=d+\mbox{\rm ad}_{u_1}$. It follows that
$d'(L_{\b,j})\in\BB^{(i+j)}$ by \eqref{equ-der-1-ebj}. By induction
on $i$, it follows that $d'$ is an inner derivation. In turn,
$d=d'-\mbox{\rm ad}_{u_1}$ is also an inner derivation.

For the other case $\a=0$, $i\neq0$, we see immediately that
$e_{0,0}=0$ by (\ref{equ-der-1-ebj}).  Applying $d$ to
$[L_{\b-1,j},L_{1,0}]= (j+(2-\b)q)L_{\b,j}$ and $[L_{\b,j},L_{-1,0}]
=-(j+(1+\b)q)L_{\b-1,j}$ respectively, we obtain two equations on
$e_{\b-1,j}$ and $e_{\b,j}$ as follows
\begin{eqnarray}
(i+j+(2-\b)q)e_{\b-1,j}+(j+q+(1-\b)(i+q))e_{1,0}&=&(j+(2-\b)q)e_{\b,j},
\nonumber\\
(i+j+(1+\b)q)e_{\b,j}+(j+q+\b(i+q))e_{-1,0}&=&(j+(1+\b)q)e_{\b-1,j}.
\nonumber
\end{eqnarray}
In particular, taking $\b=j=0$ in the first equation, we see that
\begin{equation}\label{equ-der-2-e-10e10}
e_{-1,0}+e_{1,0}=0.
\end{equation}
In turn, solving the two equations by canceling the term
$e_{\b-1,j}$, together with (\ref{equ-der-2-e-10e10}), we obtain
$i(i+2j+3q)(e_{\b,j}-\b e_{1,0})=0$, which implies that, for
$i\ne0$,
\begin{equation}\label{equ-der-2-ebj}
e_{\b,j}=\b e_{1,0}\ \ \mbox{for}\ \ \b\in\Z, j\in\Z_+.
\end{equation}
Set $u_2=(i+q)^{-1}e_{1,0}L_{0,i}\in\BB$ and let $d''=d-\mbox{\rm
ad}_{u_2}$. Then $d''(L_{\b,j})\in\BB^{(i+j)}$ by
\eqref{equ-der-2-ebj}. As in the first case, by induction on $i$, we
know that $d''$ is an inner derivation, and then $d$ is also an
inner derivation.\qed

\begin{lemm}\label{lemma-der-2st}
If $\a=i=0$, then $d$ in {\rm(\ref{equ-der-suppose})} can be written
as $d=\mbox{\rm ad}_u+\l D_0$ for some $u\in\BB$ and $\l\in\C$.
\end{lemm}

\ni{\it Proof.} \ \ Applying $d$ to
$[L_{\b,j},L_{\g,k}]=\left(\g(j+q)-\b(k+q)\right)L_{\b+\g,j+k}$, we
have
\begin{equation}\label{equ-der-3-general}
(\g(j+q)-\b(k+q))(e_{\b,j}+e_{\g,k}-e_{\b+\g,j+k})=0.
\end{equation}
Replacing $(\b,\g,k)$ by $(\b-1,1,0)$ and $(\b,-1,0)$ respectively
in \eqref{equ-der-3-general}, we have
\begin{eqnarray}
(j+(2-\b)q)(e_{\b-1,j}+e_{1,0}-e_{\b,j})&=&0,\label{equ-der-3-L10}\\
(j+(1+\b)q)(e_{\b-1,j}-e_{-1,0}-e_{\b,j})&=&0.\label{equ-der-3-L-10}
\end{eqnarray}
If $\b\neq2+jq^{-1}$, then $e_{\b,j}=e_{\b-1,j}+e_{1,0}$ by
(\ref{equ-der-3-L10}). By induction on $\b$, it follows that
\begin{eqnarray}\label{equ-der-3-ebj'}
 e_{\b,j}=
\left\{
\begin{aligned}
 &e_{\b_0-1,j}+(\b-\b_0+1)e_{1,0} && \mbox{if} \ \ \b\leq\b_0-1,\\
 &e_{\b_0,j}+(\b-\b_0)e_{1,0} && \mbox{if} \ \  \b\geq\b_0+1,
\end{aligned} \right.
\end{eqnarray}
where $\b_0:=[2+jq^{-1}]$ (the integral part of $2+jq^{-1}$). If
$\b=2+jq^{-1}=\b_0$, then
$e_{\b_0-1,j}=e_{\b_0,j}+e_{-1,0}=e_{\b_0,j}-e_{1,0}$ by
(\ref{equ-der-3-L-10}) and (\ref{equ-der-2-e-10e10}) respectively.
This, together with (\ref{equ-der-3-ebj'}), gives
$e_{\b,j}=e_{\b_0,j}+(\b-\b_0)e_{1,0}$. In particular,
$e_{0,j}=e_{\b_0,j}-\b_0 e_{1,0}$, which implies
\begin{equation}\label{equ-der-3-ebj''}
e_{\b,j}=\b e_{1,0}+e_{0,j} \ \ \mbox{for} \ \ \b\in\Z, j\in\Z_+.
\end{equation}
Furthermore, substituting \eqref{equ-der-3-ebj''} in
\eqref{equ-der-3-general} gives
$$
(\g(j+q)-\b(k+q))(e_{0,j}+e_{0,k}-e_{0,j+k})=0.
$$
Then $e_{0,j+k}=e_{0,j}+e_{0,k}$ by arbitrariness of $\b$ or $\g$.
By induction on $j$, one can derive that $e_{0,j}=je_{0,1}$, which,
together with \eqref{equ-der-3-ebj''}, gives
\begin{equation}\label{equ-der-3-ebj}
e_{\b,j}=\b e_{1,0}+j e_{0,1} \ \ \mbox{for} \ \ \b\in\Z, j\in\Z_+.
\end{equation}
Set $u_3=q^{-1}e_{1,0}L_{0,0}\in\BB$ and let $\bar{d}=d-\mbox{\rm
ad}_{u_3}-e_{0,1}D_0$, where $D_0$ is defined by
\eqref{equ-der-outer}. One will see that
$\bar{d}(L_{\b,j})\in\BB^{(j)}$ by \eqref{equ-der-3-ebj}. By Lemma
\ref{lemma-der-1st}, $\bar{d}$ is an inner derivation, and then
$d=\mbox{\rm ad}_u+e_{0,1}D_0$ for some $u\in\BB$. This completes
the proof.\qed

\begin{lemm}\label{lemma-der-3st-existence}
For every derivation $d$ of
degree $\a$, there exists some $i\in\Z_+$ such that
\eqref{equ-der-suppose} holds.
\end{lemm}
\ni{\it Proof.}\ \ Obviously, we can choose some
$i\in\Z_+$ such that
$d(L_{\a,0})\in\BB^{[i]},\,d(L_{0,1})\in\BB^{[i+1]}$ for
$\a=\pm1,\pm2$. Since $\BB$ is generated by
$\{L_{\a,0},L_{0,1}\,|\,\a=\pm1,\pm2\}$, and a derivation is
uniquely determined by its action on the generators, we see that
\eqref{equ-der-suppose} holds by induction on $j$.\qed

\begin{lemm}\label{lemma-der-4st}
For every derivation $d$, there exist derivations $d_\a$ of degree
$\a$ for all $a\in\Z$ such that $d=\sum_{\a\in\Z}d_\a$, and
furthermore, $d_\a=0$ for all but a finite number of $\a$'s.
\end{lemm}
\ni{\it Proof.}\ \ For any $\a\in\Z$, we define $d_\a$ as
follows: Let $x\in\BB_\b$ be any homogenous element of degree
$\b\in\Z$. Suppose $d(x)=\sum_{\g\in\Z}y_\g$ with $y_\g\in\BB_\g$.
Then we set $d_\a(x)=y_{\a+\b}$. This uniquely defines a linear map
$d_\a$ which can be easily verified to be a derivation of degree
$\a$. From the definition, we have
\begin{equation}\label{equ-der-deg}
d=\sum_{\a\in\Z}d_\a,
\end{equation}
 which holds in the sense that
for any $x\in\BB$, we have $d_\a(x)=0$ for all but a finite many of
$\a$'s, and $d(x)=\sum_{\a\in\Z}d_\a(x)$ (such a sum in
\eqref{equ-der-deg} is {\it summable} in this sense). This above
three lemmas show that exist $u_\a\in\BB_\a$ and some $\l_0\in\C$
such that $d_\a={\rm ad}_{u_\a}+\d_{\a,0}\l_0D_0$.

Applying \eqref{equ-der-deg} to $L_{0,0}$, by \eqref{B-block}, we
obtain that $d(L_{0,0})=-\sum_{a\in\Z}\a q u_\a$, which in
particular implies that $u_\a=0$ for all but a finite number of
$\a$'s.\qed

\vskip15pt \ni {\bf 4. \ Central extensions of
$\BB$}\setcounter{section}{4}\setcounter{theo}{0}\setcounter{equation}{0}

\ni Central extensions play an important role in the theory of Lie
algebras, since one can construct many infinite dimensional Lie
algebras by central extension and further describe the structures or
representations of these Lie algebras. The theory of universal
central extensions of Lie algebras over fields is mainly due to
Garland [\ref{G}], in which he constructs a model by using
2-cocycles. On the other hand, the cohomology groups are closely
related to the structures of Lie algebras, hence the computation of
cohomology groups seems to be important as well.

A \emph{$2$-cocycle} on $\BB$ is a $\C$-bilinear form $\psi:
\BB\times\BB\rightarrow\C$ satisfying the following conditions:
\begin{eqnarray}
{\rm (i)} && \psi(x,y)=-\psi(x,y), \nonumber\\
{\rm (ii)} && \psi([x,y],z)+\psi([y,z],x)+\psi([z,x],y)=0 \nonumber
\end{eqnarray}
for $x,y,z\in\BB$. The set of all 2-cocycles on $\BB$ is a vector
space, denoted by $Z^2(\BB,\C)$. For any $\C$-linear functions $f$
from $\BB$ to $\C$, define a 2-cocycle $\psi_f$ as follows
\begin{equation}\label{equ-2co-triv}
\psi_f(x,y)=f([x,y]),
\end{equation}
for $x,y\in\BB$, which is usually called a \emph{$2$-coboundary}, or
a \emph{trivial $2$-cocycle} on $\BB$. The set of all 2-coboundaries
is a subspace of $Z^2(\BB,\C)$, denoted by $B^2(\BB,\C)$.  We say
that two 2-cocycles $\phi, \psi$ are \emph{equivalent} if
$\phi-\psi$ is trivial. The quotient space
$$
H^2(\BB,\C)=Z^2(\BB,\C)/B^2(\BB,\C)
$$
is called the \emph{second cohomology group of $\BB$ with
coefficients in $\C$}.

As pointed in [\ref{B2}, \ref{GF}], $\Vir$ has the unique nontrivial
one-dimensional central extension, namely ${\rm
dim\ssc\,}H^2(\Vir,\C)=1$. In this section, we shall determine the
central extension of $\BB$ and the second cohomology group
$H^2(\BB,\C)$. In fact, we obtain the following analogous results of
$\BB$.

\begin{theo}\label{thm-2cocycle}
The unique nontrivial central extension of $\BB$ is given by
\begin{equation}\label{equ-2co-central-extension}
[L_{\a,i},L_{\b,j}]=\left(\b(i+q)-\a(j+q)\right)L_{\a+\b,i+j}
+\phi(L_{\a,i},L_{\b,j})c
\end{equation}
for $\a,\b\in\Z$, $i,j\in\Z_+$, where $c$ is a central element and
$\phi$ is the following non-trivial $2$-cocycle
\begin{equation}\label{equ-2co}
\phi(L_{\a,i},L_{\b,j})=\d_{\a+\b,0}\d_{i,0}\d_{j,0}\frac{\a^3-\a}{12}.
\end{equation}
Hence the second cohomology group of $\BB$ is $H^2(\BB,\C)=\C\phi$.
\end{theo}
\ni{\it Proof.} \ \ We shall prove the above theorem in a series of
lemmas. In particular, Lemmas \ref{lemma-2co-3st} and
\ref{lemma-2co-4st} imply that a 2-cocycle $\phi$ to be defined
later takes the required form (\ref{equ-2co}). This non-trivial
2-cocycle $\phi$ induces the central extension of $\BB$ as
(\ref{equ-2co-central-extension}) by taking
$c=2\phi(L_{2,0},L_{-2,0})$ (see Lemma \ref{lemma-2co-4st}(2)).\qed

First, we have
$$
[L_{\a,i},L_{0,0}]=-\a q L_{\a,i}, \ \
[L_{-1,i},L_{1,0}]=(i+2q)L_{0,i}.
$$
Let $\psi$ be any 2-cocycle. Define a linear function on $\BB$ as
follows:
\begin{equation} \label{equ-2co-def}
 f(L_{\a,i})=
\left\{
\begin{aligned}
 &-(\a q)^{-1}\psi(L_{\a,i},L_{0,0})&& \mbox{if} \ \ \a\neq0,\\
 &(i+2q)^{-1}\psi(L_{-1,i},L_{1,0})&& \mbox{otherwise.}
\end{aligned}
\right.
\end{equation}
Then $\phi=\psi-\psi_f$ is a 2-cocycle of $\BB$, which is equivalent
to $\psi$, where $\psi_f$ is the trivial 2-cocycle induced by $f$ as
in (\ref{equ-2co-triv}). Thus, by (\ref{equ-2co-def}), we
immediately have
\begin{align}
&\phi(L_{\a,i},L_{0,0}) = 0 \mbox{ \ if \ } \a\neq0,\label{equ-2co-L00,ai}\\
&\phi(L_{-1,i},L_{1,0}) = 0. \label{equ-2co-L10,-1i}
\end{align}

\begin{lemm}\label{lemma-2co-2st}
If $i\ne0$, then $\phi(L_{-2,i},L_{2,0})=0$.
\end{lemm}

\ni{\it Proof.} \ \ Applying $\phi$ to triple
$(L_{0,0},L_{1,0},L_{-1,i})$, by (\ref{equ-2co-L10,-1i}), we have
\begin{equation*}
\begin{split}
0&=
 \frac{1}{i+2q}\big(\phi([L_{-1,i},L_{0,0}],
 L_{1,0})+\phi([L_{0,0},L_{1,0}],L_{-1,i})\big) \\
&= \frac{1}{i+2q}\phi(L_{0,0},[L_{1,0},L_{-1,i}])=
 \phi(L_{0,i},L_{0,0}).
\end{split}
\end{equation*}
In turn, applying $\phi$ to triple $(L_{0,i},L_{1,0},L_{-1,0})$, it
follows from the above formula and (\ref{equ-2co-L10,-1i}) that
\begin{equation}\label{equ-2co-L-10,1i}
\begin{split}
0&=\frac{1}{i+q}\phi(L_{0,i},[L_{1,0},L_{-1,0}])\\
 &=\frac{1}{i+q}\big(\phi([L_{0,i},L_{1,0}],L_{-1,0})+
   \phi([L_{-1,0},L_{0,i}],L_{1,0})\big) \\
 &= \phi(L_{1,i},L_{-1,0}).
\end{split}
\end{equation}
So, applying $\phi$ to triple $(L_{2,0},L_{-1,i},L_{-1,0})$, by
(\ref{equ-2co-L10,-1i}) and (\ref{equ-2co-L-10,1i}), we have
\begin{equation*}
\begin{split}
0&=\phi([L_{-1,0},L_{2,0}],L_{-1,i})+
   \phi([L_{2,0},L_{-1,i}],L_{-1,0}) \\
 &=\phi(L_{2,0},[L_{-1,i},L_{-1,0}])=i\phi(L_{-2,i},L_{2,0}),
\end{split}
\end{equation*}
which gives the desired conclusion.\qed

\begin{lemm}\label{lemma-2co-3st}
For $\a, \b\in\Z$, $i, j\in\Z_+$, we have
\baselineskip3pt\lineskip7pt\parskip-3pt
\begin{itemize}\parskip-3pt
  \item[{\rm(1)}]\ $\phi(L_{\a,i},L_{\b,j})=0$ if $\a+\b\neq0$;
  \item[{\rm(2)}]\ $\phi(L_{\a,i},L_{-\a,j})=0$ if $i\neq j$.
\end{itemize}
\end{lemm}

\ni{\it Proof.} \ \ If $\a+\b\neq0$, applying $\phi$ to triple
$(L_{0,0},L_{\a,i},L_{\b,j})$, by (\ref{equ-2co-L00,ai}), we see
\begin{equation*}
\begin{split}
0&=\frac{1}{(\a+\b)q}\phi(L_{0,0},[L_{\a,i},L_{\b,j}])\\
 &=\frac{1}{(\a+\b)q}\big(\phi([L_{0,0},L_{\a,i}],L_{\b,j})
+\phi([L_{\b,j},L_{0,0}],L_{\a,i})\big)\\
 &=\phi(L_{\a,i},L_{\b,j}),
\end{split}
\end{equation*}
which gives part (1). Applying $\phi$ to triples
$(L_{1,0},L_{\a,i},L_{-1-\a,j})$ and
$(L_{-1,0},L_{1+\a,i},L_{-\a,j})$ respectively, by
(\ref{equ-2co-L10,-1i}) and (\ref{equ-2co-L-10,1i}), we obtain
following two equations
\begin{equation*}\label{equ-2co-fromL10}
\begin{split}
0&=\phi(L_{1,0},[L_{\a,i},L_{-1-\a,j}])\\
 &=\phi([L_{1,0},L_{\a,i}],L_{-1-\a,j})+
   \phi([L_{-1-\a,j},L_{1,0}],L_{\a,i})\\
 &=((\a-1)q-i)\phi(L_{1+\a,i},L_{-1-\a,j})-
   ((\a+2)q+j)\phi(L_{\a,i},L_{-\a,j}),
\end{split}
\end{equation*}
\begin{equation*}\label{equ-2co-fromL-10}
\begin{split}
0&=\phi(L_{-1,0},[L_{1+\a,i},L_{-\a,j}])\\
 &=\phi([L_{-1,0},L_{1+\a,i}],L_{-\a,j})+
   \phi([L_{-\a,j},L_{-1,0}],L_{1+\a,i})\\
 &=((\a+2)q+i)\phi(L_{\a,i},L_{-\a,j})-
   ((\a-1)q-j)\phi(L_{1+\a,i},L_{-1-\a,j}).
\end{split}
\end{equation*}
Multiplying the first equation by $(\a-1)q-j$, the second one by
$(\a-1)q-i$, and then adding both results together, we deduce
$$
(i-j)(i+j+3q)\phi(L_{\a,i},L_{-\a,j})=0,
$$
which immediately gives part (2) since $i,j\in\Z_+$ and $q>0$.\qed

\begin{lemm}\label{lemma-2co-4st}
For $\a\in\Z$, $i\in\Z_+$, we have
\baselineskip3pt\lineskip7pt\parskip-3pt
\begin{itemize}\parskip-3pt
  \item[{\rm(1)}]\ $\phi(L_{\a,i},L_{-\a,i})=0$ if $i\ne0$;
  \item[{\rm(2)}]\ $\phi(L_{\a,0},L_{-\a,0})=\frac{\a^3-\a}{6}\phi(L_{2,0},L_{-2,0})$.
\end{itemize}
\end{lemm}

\ni{\it Proof.} \ \ Applying $\phi$ to triple
$(L_{2,0},L_{-2-\a,i},L_{\a,i})$, we get
\begin{equation}\label{equ-2co-fromL20}
\begin{split}
&2(i+q)(\a+1)\phi(L_{2,0},L_{-2,2i})=\phi(L_{2,0},[L_{-2-\a,i},L_{\a,i}])\\
&=\phi([L_{2,0},L_{-2-\a,i}],L_{\a,i})+
  \phi([L_{\a,i},L_{2,0}],L_{-2-\a,i})\\
&=((\a+4)q+2i)\phi(L_{\a,i},L_{-\a,i})-
  ((\a-2)q-2i)\phi(L_{2+\a,i},L_{-2-\a,i}).
\end{split}
\end{equation}
Furthermore, applying $\phi$ to triples $(L_{1,0}, L_{\a,i},
L_{-1-\a,i})$ and $(L_{1,0}, L_{1+\a,i}, L_{-2-\a,i})$ respectively
(here we omit the details), by (\ref{equ-2co-L10,-1i}), we have
\begin{eqnarray}
((\a+2)q+i)\phi(L_{\a,i},L_{-\a,i})
&=&((\a-1)q-i)\phi(L_{1+\a,i},L_{-1-\a,i}),\label{equ-2co-fromL10'} \\
((\a+3)q+i)\phi(L_{1+\a,i},L_{-1-\a,i})&=&(\a
q-i)\phi(L_{2+\a,i},L_{-2-\a,i}). \label{equ-2co-fromL10''}
\end{eqnarray}
Solving the three equations
\eqref{equ-2co-fromL20}--\eqref{equ-2co-fromL10''} by canceling the
common terms $\phi(L_{1+\a,i},L_{-1-\a,i})$ and
$\phi(L_{2+\a,i},L_{-2-\a,i})$, we obtain
$$
(i+2q)(2i+3q)\phi(L_{\a,i},L_{-\a,i}) =(\a+1)(\a
q-i)((\a-1)q-i)\phi(L_{2,0},L_{-2,2i}).
$$
If $i\ne0$, then the above formula gives part (1) by Lemma
\ref{lemma-2co-2st}. Otherwise, $i=0$ and part (2) clearly
holds.\qed\vskip7pt

\small \ni{\bf Acknowledgement} The authors would like to thank
Professor Yucai Su for instructions and helps.

\small\def\bf{}
\parskip=-1pt\baselineskip=3pt\lineskip=3pt

\end{document}